\documentclass[11pt]{amsart}
\numberwithin{equation}{section}

\title[Equivalence of probability spaces]
{On the equivalence of probability spaces}
\usepackage{amsthm}
\usepackage{amsmath,mathrsfs}
\usepackage{color}
\usepackage{amssymb}
\usepackage{hyperref}
\usepackage{enumerate}
\usepackage{latexsym,array}
\usepackage{amsfonts}
\usepackage{shadow}

\newtheorem{Pa}{Paper}[section]
\newtheorem{Tm}[Pa]{{\bf Theorem}}
\newtheorem{La}[Pa]{{\bf Lemma}}
\newtheorem{Dn}[Pa]{{\bf Definition}}
\newtheorem{Cy}[Pa]{{\bf Corollary}}

\newtheorem{Rk}[Pa]{{\bf Remark}}
\newtheorem{Pn}[Pa]{{\bf Proposition}}

\newtheorem{Ex}[Pa]{{\bf Example}}

\newtheorem{Nn}[Pa]{{\bf Notation}}

\newcommand{\w}{\omega}
\newcommand{\e}{\varepsilon}

\author[D. Alpay]{Daniel Alpay}
\address{(DA) Department of Mathematics\\
Ben-Gurion University of the Negev\\
Beer-Sheva 84105, Israel} \email{dany@math.bgu.ac.il}

\author[P. Jorgensen]{Palle Jorgensen}
\address{(PJ)  Department of Mathematics\\
 University of Iowa.
Iowa City, IA 52242 USA}
\email{palle-jorgensen@uiowa.edu}

\author[D. Levanony]{David Levanony}
\address{(DL) Department of Electrical Engineering\\ Ben Gurion
University of the Negev,\\ Beer Sheva 84105, Israel}
 \email{levanony@ee.bgu.ac.il}

\def\R{\mathbb R}

\def\(s){\mathscr S(\R\times\R)}


\subjclass[2010]{}

\thanks{D. Alpay and P. Jorgensen thank the Binational
Science Foundation Grant number 2010117. One of the authors (PJ)
thanks colleagues at Ben-Gurion University for kind hospitality,
and for many very fruitful discussions. Part of this work was
done while PJ visited BGU in May and June 2014. D. Alpay thanks
the Earl Katz family for endowing the chair which supported his
research.}
\begin{document}
\parindent 0cm
\begin{abstract}
For a general class of Gaussian processes  $W$, indexed by a sigma-algebra $\mathscr F$  of a general
measure space
$(M,\mathscr F, \sigma)$, we give necessary and sufficient conditions for the validity of a quadratic variation 
representation for such Gaussian processes, thus recovering $\sigma(A)$, for $A\in\mathscr F$, as a quadratic 
variation of $W$ over $A$. 
We further provide a harmonic analysis representation for this general class of processes.
We apply these two results to: $(i)$ a computation of generalized Ito-integrals; and $(ii)$ a
proof of an explicit, and measure-theoretic equivalence formula, realizing an equivalence between the two approaches to 
Gaussian processes, one where the choice of sample space is the traditional path-space, and the other 
where it is Schwartz' space of tempered distributions.
\end{abstract}
\maketitle
\tableofcontents
\section{Introduction}
\setcounter{equation}{0}
Any stochastic process must be realized in a probability space, a triple made up of a sample space  $\Omega$,
a sigma-algebra
$\mathscr B$, and a choice of probability measure $P$ defined on
$(\Omega, \mathscr B)$. In research papers so far, covering a
general family of stationary-increment processes, the
choice of sample space  $\Omega$ has typically been   $\mathcal
S^\prime$ (the space of tempered distributions in the sense of
Laurent Schwartz). However for applications, a better choice for
$\Omega$ (the sample space) is clearly the continuous functions,
$C(\mathbb R)$  if the process is indexed by time
(sample paths for the process.)  In a setting more general than
$\mathcal S^\prime$ versus $C(\mathbb R)$  we present an
explicit measure-isomorphism between the two choices. It will be
given by an explicit formula, will be well suited for
computations; and our formulas are new even in the case of
$\mathcal S^\prime$ versus $C(\mathbb R)$. In addition to this
isomorphism, we offer a number of applications to stochastic
integration.\\

In our discussion below and in our proofs, we will be making use
of the following earlier papers
\cite{MR0622034,MR0370681,MR2793121}.   The literature on
families of Gaussian processes is vast, both in pure and in
applied areas each with a different viewpoint. Our present motivation stems
from a certain harmonic analysis, from spectral theoretic
computations, and from applications in quantum theory. In our
analysis, we study a particular family of Gaussian processes
(see sections \ref{sec2} and \ref{sec4}) and, in our proofs, we
utilize of tools from a host of areas. Specifically, we utilize families of
operators in Hilbert spaces, with emphasis on explicit reproducing
kernel Hilbert spaces. To refer the reader to useful
references, we mention the following
papers/books whose theme is close to that adopted here:
\cite{ajnfao,MR2793121,MR2788708,MR2444857,MR2397796,MR1810995,MR2805533,MR2459226},
and for papers/books stressing the Hilbert space, see e.g.,
\cite{MR2521338,MR2105321,MR2918388,MR2797301,MR2766605,MR2683408,MR3081915,MR2920083,
MR2662722,MR2480780,MR2518590}, and the papers cited therein. We make use of basic tools from harmonic analysis, and Gaussian
processes in Hilbert space; for background references, see e.g.,
\cite{MR1216518,lifs,MR1943282,MR2932674,MR2728698,MR2512800}.\\

The paper is organized as follows: In the first half of the paper
(sections \ref{sec2}-\ref{sec4}), we introduce a general class of
Gaussian processes indexed by a sigma-algebra of a general but fixed
measure space. This material will aid us in two ways, in the second half of the paper,
sections \ref{sec6}-\ref{sec7}.  First it will unify our approach to
generalized Ito-integrals, subsequently studied. Secondly, our general theory will be
used in the proofs of our results covering our main theme, i.e., setting up a
measure-theoretic equivalence between the two approaches to the formulation of Gaussian
processes outlined above. This thread culminates with Theorem \ref{thm82}.

\section{Preliminaries}
\setcounter{equation}{0}
\label{sec2}
 Below we briefly sketch the framework for the particular class of Gaussian processes to be studied in detail.
We first gather some notations and definitions needed in the sequel.
We already mentioned the Schwartz space $\mathcal S$ and its dual $\mathcal S^\prime$. The duality between the two spaces
will be denoted as
\begin{equation}
\label{pairing}
\langle \xi,\varphi\rangle,\quad \varphi\in\mathcal S\quad\text{\and}\quad \xi\in\mathcal S^\prime,\,\,\text{the distribution $\xi$ applied to
$\varphi$.}
\end{equation}
\begin{Dn}
\label{def21}
Let $M=\mathbb R$ and let $\mathscr F=\mathscr
B(\mathbb R)$ denote the Borel sigma-algebra. A measure $\sigma$ on $\mathscr B(\mathbb R)$ is called
tempered if it satisfies
\begin{equation}
\label{8} \int_{\mathbb R}\frac{d\sigma(u)}{(u^2+1)^p}<\infty
\end{equation}
for some $p\in\mathbb N_0$.\\
The measure $\sigma$ will be called symmetric if
\begin{equation}
\label{neww2}
\sigma(A)=\sigma(-A),\quad\forall A\in\mathscr B.
\end{equation}
\end{Dn}

\begin{Dn} A cylinder set in $\mathcal S^\prime$ is a subset of $\mathcal S^\prime$ of the form
\begin{equation}\label{eqdefcyl}
\left\{\xi\in\mathcal S^\prime\,;\, (\langle \xi,\varphi_1\rangle,\ldots, \langle \xi,\varphi_n\rangle)\in A\right\}
\end{equation}
where $\varphi_1,\ldots,\varphi_n$ are preassigned in the Schwartz space $\mathcal S$, and  $A\subset \mathbb R^n$ is an 
open subset of $\mathbb R^n$.
\label{defcyl}
\end{Dn}

We denote the Fourier transform by
\begin{equation}
\widehat{f}(\lambda)=\int_{\mathbb R}e^{-i\lambda x}f(x)dx,\quad \text{for}\quad f\in\mathbf L_2(\mathbb R).
\end{equation}

\begin{Tm}
\label{tm32}
Let $\sigma$ be fixed. There is a real-valued Gaussian process
$X^{(\sigma)}$ indexed by the space of Schwartz functions
$\mathcal S$ and realized on a probability space $(\mathcal
S^\prime,\mathscr C,P^{(\sigma)})$, where $\mathscr C$ denotes the
cylinder algebra, such that, with $\mathbb
E^{(\sigma)}[u]=\int_{\mathcal S^\prime}udP^{(\sigma)}$ and with the notation \eqref{pairing},
\[
X_\varphi^{(\sigma)}(\xi)=\langle\xi,\varphi\rangle,
\]
we have\smallskip

$(1)$
\[
\mathbb
E^{(\sigma)}\left[X^{(\sigma)}_\varphi\right]=0,\quad\forall
\varphi\in\mathcal S.
\]

and\smallskip

$(2)$
\begin{equation}\label{090714}
\mathbb
E^{(\sigma)}\left[e^{iX^{(\sigma)}_\varphi}\right]=e^{-\frac{1}{2}\int_{\mathbb
R}|\widehat{\varphi}(u)|^2d\sigma(u)}.
\end{equation}
\end{Tm}

We will return to a detailed construction    of this process in Sections
\ref{sec4}-\ref{sec7}, but we first need some technical
preparation.
\begin{proof}[Proof of Theorem \ref{tm32}] In the construction of the Gaussian process
\begin{equation}
\label{gauss090714}
(X_\varphi^{(\sigma)})_{\varphi\in\mathcal S}\quad\text{on}\quad (\mathcal S^\prime,\mathscr C, P^{(\sigma)})
\end{equation}
in Theorem \ref{tm32}, we "invert" a transform. Specifically, we
consider the right hand side in \eqref{090714} as a continuous
and positive definite function on $\mathcal S$, and we apply
Minlos' theorem to verify the existence of the Gaussian process
\eqref{gauss090714}.
\end{proof}
\begin{Rk}
\label{rknew1}
{\rm By this approach, we do not get
explicit formulas for \eqref{gauss090714} that can be computed on
sample-paths. This situation is remedied in Theorems \ref{tm61}
and \ref{thm82} below; see especially formulas
\eqref{new17}--\eqref{new19}.}
\end{Rk}
\section{ Gaussian processes}
\setcounter{equation}{0} \label{sec3}

This section is divided into four parts. In the first, our starting point is a given sigma-finite measure space $(M,\mathscr F, \sigma)$, and we outline properties of a Gaussian process  $W^{(\sigma)}$  which is indexed canonically by the sets from the sigma-algebra of $M$. In subsection \ref{sec3.2},
we prove an inversion-formula: We show that if $\sigma$ is refinable (Definition \ref{newdf1}), then this measure $\sigma$  may be recovered as a quadratic variation computed from $W^{(\sigma)}$ . Moreover, with this, we formulate and prove a generalized Ito-lemma for $W^{(\sigma)}$. And, as a result, we obtain a general stochastic integration. 
We use this in the subsequent Section 
\ref{sec3.3} to introduce a coordinate system (Theorem \ref{thm9}) on the probability space  $\Omega$ carrying the process. 
Moreover, we obtain (Corollary \ref{cor12}) a generalized Fourier transform in $\mathbf L_2(\Omega, P)$. We write the
proofs in the case of real valued functions, and real valued
random variables, but the arguments apply {\it mutatis mutandis}
for the complex case as well.

\subsection{Measure space associated Gaussian processes}

\begin{Dn}
\label{def1} Let $(M,\mathscr F,\sigma)$ be a sigma-finite measure
space, meaning that there exists an increasing sequence
$A_1\subset A_2\subset\cdots$ of elements of $\mathscr F$, such
that $\cup_{n=1}^\infty A_n=M$ and $\sigma(A_n)<\infty$ for all
$n\in\mathbb N$.
We say that the Gaussian process $W=(W_A^{(\sigma)})_{A\in\mathscr F}$ is {\rm associated with $(M,\mathscr F, \sigma)$} if there is a probability space $(\Omega,\mathcal B,P)$ such that the following two conditions hold:\\
$(i)$ Setting $\mathbb E[X]=\int_\Omega X(\w)dP(\w)$, we have
\begin{eqnarray}
\label{1}
\mathbb E[W_A]&=&0,\quad\forall A\in\mathscr F,\\
\label{2}
\mathbb E[W_AW_B]&=&\sigma(A\cap B),\quad \forall A,B\in\mathscr F.
\end{eqnarray}
$(ii)$ For every finite sequence $A_1,A_2,\ldots,A_n$ of elements
of $\mathscr F$, the collection of random variables
\[
W_{A_1},W_{A_2}, \ldots,W_{A_n}
\]
is jointly Gaussian with
covariance
\begin{equation}
\label{3}
\mathbb E[W_{A_i}W_{A_j}]=\sigma(A_i\cap A_j),\quad i,j=1,\ldots, n.
\end{equation}
\end{Dn}

\begin{La}
For every sigma-finite measure space $(M,\mathscr F,\sigma)$, an
associated Gaussian process exists. \label{lem1}
\end{La}

\begin{proof}The function $K(A,B)=\sigma(A\cap B)$ is positive definite on $\mathscr F$;  and so it is
the covariance function of a zero-mean Gaussian process with covariance function $\sigma(A\cap B)$,
as follows from \cite[pp. 466-467]{MR58:31324b}.
\end{proof}

We note that in the present setting,  reference \cite{MR0622034} is also helpful for a general discussion of positive definite kernels and Gaussian processes.

\begin{Ex} {\rm The Cameron-Martin space revisited:}
Let $(M,\mathscr F,\sigma)=(\mathbb R_+,\mathscr B, \lambda)$,
where $\mathbb R_+$ is the half-line $[0,\infty)$, $\mathscr B$
denotes the Borel subsets, and $\lambda$ is the restriction of the
Lebesgue measure to $[0,\infty)$. Then the reproducing kernel
Hilbert space associated with the positive definite function
\[
K(A,B)=\lambda(A\cap B),\quad A,B\in\mathscr B,
\]
is the Cameron-Martin Hilbert space $\mathcal H_1^{(\lambda)}$ where
\[
\mathcal H_1^{(\lambda)}=\left\{f\in\mathbf L_2(\mathbb R_+)\,|\, f^\prime\in \mathbf L_2(\mathbb R_+),\,\, and\,\, f(0)=0\right\}
\]
and
\[
\|f\|^2_{\mathcal H_1^{(\lambda)}}=\int_0^1|f^\prime(x)|^2dx,\quad f\in \mathcal H_1^{(\lambda)}.
\]
\label{cameron}
\end{Ex}
\begin{proof}[Discussion of Example \ref{cameron}:]
For $s,t\in[0,\infty)$ set $A=[0,s]$ and $B=[0,t]$. Then
\[
\lambda (A\cap B)=s\wedge t=k(s,t)=k_s(t),
\]
and
\[
\langle f,k_s\rangle_{\mathcal H_1^{(\lambda)}}=f(s),\quad \forall f\in \mathcal H_1^{(\lambda)}\quad \text{and}\quad s\in[0,\infty)
\]
\end{proof}

\begin{Rk}
\label{rkvu160714}
{\rm The conclusion in the example above is still valid if the measure $\sigma$ for $(\mathbb R_+,\mathscr B,\sigma)$ is 
tempered; see Definition \ref{def21}. The modifications are as follows:\\
$(i)$ With $\sigma$ given, we get $W^{(\sigma)}$ and its covariance kernel. The corresponding reproducing kernel Hilbert 
space $RK(\sigma)$ is described as follows:
\[
RK(\sigma)=\left\{F\,\,:\,\, F(t)=\int_0^tf(x)d\sigma(x)\,\,\text{for some $f\in\mathbf L_2(\sigma)$}\right\}
\]
with norm
\[
\|F\|^2_{RK(\sigma)}=\int_0^\infty |f(x)|^2d\sigma(x),
\]
and\\
$(ii)$ the reproducing kernel for $RK(\sigma)$ is
\[
k^{(\sigma)}(t,s)=\sigma([0,t\wedge s]).
\]
}
\end{Rk}
\begin{La}
\label{lemma3} If $(M,\mathscr F,\sigma)$ is a sigma-finite
measure space, and $(W_A^{(\sigma)})_{A\in\mathscr F}$ is an
associated Gaussian process on the probability space
$(\Omega,\mathcal B,P)$, then the Ito integral
\[
\int_Mf(x)dW^{(\sigma)}(x)\,\,\in\,\,\mathbf L_2(\Omega,\mathcal
B,P)
\]
exits for all $f\in\mathbf L_2(M,\mathscr F,\sigma)$, and the isometry
\begin{equation}
\label{4} \mathbf
E\left[\big|\int_Mf(x)dW^{(\sigma)}(x)\big|^2\right]=\int_M|f(x)|^2d\sigma(x)
\end{equation}
holds.
\end{La}

\begin{proof}
Let $f_s$ be a generic simple function, that is
\begin{equation}\label{5}
f_s(x)=\sum_{k=1}^nc_k\chi_{A_k}(x),\quad n\in\mathbb N,\quad A_1,\ldots,A_n\in\mathscr F\quad\text{and}\quad c_1,\ldots,
c_n\in\mathbb R,
\end{equation}
where $\chi_{A_k}$ denotes the appropriate indicator function with
the sets $A_1,\ldots,A_n$ being such that $A_j\cap A_k=\emptyset$ for
$j\not= k$, and $c_1,\ldots,c_n\in\mathbb R$. Set
\begin{equation}
\label{6}
\int_Mf_s(x)dW^{(\sigma)}(x)=\sum_{k=1}^nc_kW_{A_k}^{(\sigma)}.
\end{equation}
We claim that \eqref{4} holds. The desired conclusion \eqref{4} follows
for all $f\in\mathbf L_2(M,\mathscr F,\sigma)$ as the simple functions
are dense in this latter space and since \eqref{4} is a densely defined isometry between Hilbert spaces, and thus has a unique everywhere
and continuous isometric extension. Indeed (and with the limit meaning approximation by simple functions),
\[
\begin{split}
\mathbb E\left[\big|\int_Mf(x)dW^{(\sigma)}(x)\big|^2\right)&=\lim\sum_{j,k=1}^nc_jc_k{\mathbb E}\left[W_{A_j}^{(\sigma)}
W_{A_k}^{(\sigma)}\right] \\
&=\lim\sum_{j,k=1}^nc_jc_k\sigma(A_j\cap A_k)\\
&=\lim\sum_{k=1}^nc_k^2\sigma(A_k)\\
&=\int_M|f(x)|^2d\sigma(x)
\end{split}
\]
as claimed.
\end{proof}

%

\begin{Rk}
\label{newrk}
{\rm A refinement of the proof of Lemma \ref{lemma3}, for deterministic functions, also yields the case of Ito-integral and Ito-isometry, for 
stochastic adapted process.
The purpose of the extension is to set the stage for the Ito lemma (Corollary \ref{cy322}) where the reasoning relies on the generality of stochastic adapted processes, indexed by a given measure space.
But for pedagogical reasons, we have chosen to first state the result in the “easier” special case of the Wiener integration for deterministic functions.  Below we now point out how the proof of this case carries over, with suitable modifications,  to the case of Ito-integration, and we establish the Ito-isometry for the case of a family of processes which we call $\mathcal F$-adapted. See Definition \ref{lalaldo} and Lemma \ref{lalala}.
}
\end{Rk}

\begin{Cy}
\label{cy3b} Let $\left(W_A^{(\sigma)}\right)_{A\in\mathscr F}$
be a Gaussian process defined on the probability space
$(\Omega,\mathscr B,P)$, and let $I=(a,b]$ be a finite interval.
Take an $A\in\mathscr F$ such that $0<\sigma(A)<\infty$. Then
\[
P\left(\left\{\w\in\Omega\,|\, W^{(\sigma)}_A(\w)\in
I\right\}\right)=P\left(\left\{a<W_A^{(\sigma)}\le
b\right\}\right)=\frac{1}{\sqrt{2\pi}}\int_{\frac{a}{\sqrt{\sigma(A)
}}}^{\frac{b}{\sqrt{\sigma(A)}}}e^{-\frac{x^2}{2}}dx.
\]
\end{Cy}

\begin{Cy}
\label{cyvu}
Let $W^{(\sigma)}(f)$, $f\in\mathbf L_2(\sigma)$ be as in 
\eqref{7}. Then,
\[
\mathbb E_{P_\sigma}\left[W^{(\sigma)}(f_1)W^{(\sigma)}(f_2)
\right]=\langle f_1,f_2\rangle_{\sigma}\left(=\int_Mf_1(u)f_2(u)
d\sigma(u)\right),\quad\forall f_1,f_2\in\mathbf L_2(\sigma).
\]
\end{Cy}
\begin{proof}
This is the polarization of
\[
\mathbb E_{P_\sigma}\left[\left(W^{(\sigma)}(f)\right)^2\right]
=\|f\|^2_\sigma.
\]
\end{proof}

\begin{Nn} (The Wiener integral)
With $(M,\mathscr F,\sigma)$ and $(W_A^{(\sigma)})_{A\in\mathscr
F}$ as above, we set for $f\in\mathbf L_2(M,\mathscr F,\sigma)$
\begin{equation}
\label{7} W^{(\sigma)}(f)=\int_MfdW^{(\sigma)}.
\end{equation}
\label{n4}
\end{Nn}

Let $(M,\mathscr F,\sigma)$ be as above with $W^{(\sigma)}$ being the associated Gaussian process. Let $H_0,H_1,\ldots$ denote the Hermite
polynomials, defined by the generating function
\begin{equation}
e^{zx-\frac{z^2}{2}}=\sum_{n=0}^\infty \frac{z^n}{n!}H_n(x).
\label{te1}
\end{equation}
Let $\psi\in\mathcal S$ be a (real-valued) Schwartz function, with Hermite expansion
\begin{equation}
\label{te2}
\psi(x)=\sum_{n=0}^\infty c_nH_n(x),\quad x\in\mathbb R,
\end{equation}
and set
\begin{equation}
\label{te3}
[\psi](x)\stackrel{\rm def.}{=}\sum_{n=0}^\infty n!c_n^2x^n,\quad x\in\mathbb R.
\end{equation}
With these notations we can now state:

\begin{Cy}
With $W^{(\sigma)}$, $\psi$ and $[\psi]$ as above we have the following: Let $f_1,f_2\not=0$ be in $\mathbf L_2(\sigma)$. Then,
\begin{equation}
\label{te4}
\mathbb E_{P_\sigma}\left[\psi\left(\frac{1}{\|f_1\|_\sigma}  W^{(\sigma)}(f_1)\right)
\left(\psi\left(\frac{1}{\|f_2\|_\sigma}  W^{(\sigma)}(f_2)\right)\right)\right]=[\psi]\left(\frac{\langle f_1,f_2\rangle_\sigma}{
\|f_1\|_\sigma\|f_2\|_\sigma}\right).
\end{equation}
\end{Cy}
\begin{proof} By Lemma \ref{lemma3} applied to the $\mathcal N(0,1)$ random variables $\frac{1}{\|f_1\|_\sigma}  W^{(\sigma)}(f_1)$ and $\frac{1}{\|f_2\|_\sigma}  W^{(\sigma)}(f_2)$ we have
\begin{equation}
\label{tset}
\mathbb E_{P_\sigma}\left[\left(H_n\left(\frac{1}{\|f_1\|_\sigma}  W^{(\sigma)}(f_1)\right)\right)
\left(H_k\left(\frac{1}{\|f_2\|_\sigma}  W^{(\sigma)}(f_2)\right)\right)\right]=n!\delta_{k,n}\left(\frac{\langle f_1,f_2\rangle_\sigma}{
\|f_1\|_\sigma\|f_2\|_\sigma}\right)^n.
\end{equation}
Consider now the left hand side of \eqref{te4}, and replace $\psi$ by its expansion \eqref{te2}. Taking into
account \eqref{tset} we get the right hand side of \eqref{te4}.
\end{proof}
In the computation of \eqref{tset}, we use the following fact about Gaussian vectors in $\mathbb R^2$.

\begin{La}
Let $c\in\mathbb R$, $|c|<1$, and let $\gamma_2^{(c)}$ be the $\mathbb R^2$-Gaussian joint density with covariance matrix 
$\begin{pmatrix}1&c\\c&1
\end{pmatrix}$, then for the Hermite functions $H_n$ and $H_k$ we have:
\begin{equation}
\iint_{\mathbb R^2}H_n(x)H_k(y)d\gamma_2^{(c)}(x,y)=\delta_{n,k}n!c^n,\quad n,k\in\mathbb N_0.
\end{equation}
\end{La}
The proof is by direct computation, making use of the Ornstein-Uhlenbeck semigroup $P_t$ given by
\begin{equation}
(P_th)(x)=\int_{\mathbb R}h(e^{-t}x+\sqrt{1-e^{-2t}}y)d\gamma_1(y),\quad h\in\mathbf L_2(\mathbb R,d\gamma_1).
\end{equation}
For details see for instance \cite[pp. 9-28]{MR2962301}.
\begin{Cy}
Let $W^{(\sigma)}$, $\psi$ and $[\psi]$ as above, and let $A,B\in\mathscr F$ be such that $0<\sigma(A)<\infty$ and
$0<\sigma(B)<\infty$. Then,
\begin{equation}
\label{te4b}
\mathbb E_{P_\sigma}\left[\left(\psi\left(\frac{1}{\sqrt{\sigma(A)}}  W^{(\sigma)}_A\right)\right)
\left(\psi\left(\frac{1}{\sqrt{\sigma(B)}}  W^{(\sigma)}_B\right)\right)\right]=[\psi]\left(\frac{\sigma(A\cap B)}{
\sqrt{\sigma(A)\sigma(B)}}\right).
\end{equation}
\end{Cy}

\begin{Dn}
Let $(M,\mathscr F,\sigma)$ and  $(W_A^{(\sigma)})_{A\in\mathscr
F}$ be as in Definition \ref{def1}. We say that a (finite or)
countable family $A_1,A_2,\ldots\in\mathscr F$ is a partition of
$A$ if $A=\cup_{k\in\mathbb N}A_k$, and $A_k\cap A_n=\emptyset$
for $k\not= n$.
\end{Dn}

\begin{Cy}
Let $\sigma$ be a tempered Borel measure on $\mathbb R$, and let $(X_\varphi^{(\sigma)})_{\varphi\in\mathcal S}$ be the corresponding Gaussian process, namely $X_\varphi^{(\sigma)}(\xi)=\langle\xi,\varphi\rangle$. Then,
\begin{equation}
\label{090714b} X_{\varphi}^{(\sigma)}=\int_{\mathbb
R}\widehat{\varphi}(u)dW_u^{(\sigma)},
\end{equation}
where $W^{(\sigma)}$ in \eqref{090714b} is the process in Definition \ref{def1}, and the integral is the generalized Ito 
integral of Lemma \ref{lemma3}.
\end{Cy}
\begin{proof}
Pick a finite partition $\left\{x_j\right\}$ of $\mathbb R$ such that the integral in \eqref{090714b} is approximated by
\[
\sum_j\widehat{\varphi}(x_j)W^{(\sigma)}_{[x_j,x_{j+1})}.
\]
The approximation is in the $L_2$ norm in $\mathbf L_2(\mathcal S^\prime,P_\sigma)$. Using
\[
\prod_j\mathbb E_{P_\sigma}\left[e^{i\widehat{\varphi}(x_j)W^{(\sigma)}_{[x_j,x_{j+1})}}\right]=e^{-\frac{1}{2}\sum_j|\widehat{\varphi}
(x_j)|^2\sigma\left([x_j,x_{j+1)}\right)},
\]
we conclude, upon passing to the limit of refinements, that
\begin{equation}
\label{newnew2}
\mathbb E_\sigma\left[e^{iX_{\varphi}^{(\sigma)}}\right]=e^{-\frac{1}{2}\int_{\mathbb R}|\widehat{\varphi}(u)|^2d\sigma(u)}.
\end{equation}
Equivalently, the Gaussian process defined in \eqref{090714b},
satisfies equation $(2)$ of Theorem \ref{tm32}. By
Minlos' theorem,  the process in $(2)$ is uniquely determined
by the characteristic function $e^{-\frac{1}{2}\int_{\mathbb
R}|\widehat{\varphi}(u)|^2d\sigma(u)}$.
\end{proof}

\subsection{Quadratic variation}
\label{sec3.2}
\begin{La}
Let $(M,\mathscr F,\sigma)$ be a measure space with $\sigma $ being sigma-finite, and let $(W_A^{(\sigma)})_{A\in\mathscr F}$
be the Gaussian process of Definition \ref{def1}. Then for its moments we have:\\
$(1)$ The odd moments vanish,
\[
E_{P_\sigma}\left[(W_A^{(\sigma)})^{2k+1}\right]=0,\quad k\in\mathbb N_0,
\]
and\\
$(2)$ the even moments are given by:
\[
E_{P_\sigma}\left[(W_A^{(\sigma)})^{2k}\right]=(2k-1)!!(\sigma(A))^k,\quad k\in\mathbb N,
\]
for all $A\in\mathscr F$ such that $\sigma(A)<\infty$ where,
\[
(2k-1)!!=\frac{(2k)!}{2^kk!}.
\]
\label{toto}
\end{La}
We note that in particular we have:
\[
\begin{split}
\mathbb E_{P_\sigma}[W_A^{(\sigma)}]&=0,\\
\mathbb E_{P_\sigma}[(W_A^{(\sigma)})^2]&=\sigma(A),\\
\mathbb E_{P_\sigma}[(W_A^{(\sigma)})^4]&=3(\sigma(A))^2.\\
\end{split}
\]
\begin{proof}[Proof of Lemma \ref{toto}] The asserted moment expressions in the lemma follow from comparing
powers in the moment generating function
\[
\mathbb E_{P_\sigma}\left[e^{itW_A^{(\sigma)}}\right]=e^{-\frac{t^2\sigma(A)}{2}},
\]
valid for all $A\in\mathscr F$ such that $\sigma(A)<\infty$.
\end{proof}

\begin{Pn}
\label{prop37} Let the process $(W_A^{(\sigma)})_{A\in\mathscr F}$
be realized on $\mathbf L_2(\Omega,\mathscr B,P)$ as outlined in
Lemma \ref{lemma3}, and let $\mathbb E$ be the expectation operator
defined by $P$. Then for any $A\in\mathscr F$ and every
partition $(A_k)_{k\in\mathbb N}$ of $A$, the following identity
holds:
\begin{equation}
\label{new2} \mathbb
E\left[|\sigma(A)-\sum_{k}\left(W_{A_k}^{(\sigma)}\right)^2|^2\right]=
2\sum_{k}(\sigma(A_k))^2=2\sum_{k}
\left(\mathbb E\left[(W_{A_k}^{(\sigma)})^2\right]\right)^2.
\end{equation}
\end{Pn}
\begin{proof}
We compute the left hand side in \eqref{new2} with the use of
Lemma \ref{lemma3} as follows: We first note that the random
variables $W^{(\sigma)}_{A_k}$ and $W^{(\sigma)}_{A_n}$ are
independent when $k\not=n$. This follows from \eqref{3} together with the
fact that they are Gaussian; see Corollary \ref{cy3b}. Set
$s=\sigma(A)$ and $s_k=\sigma(A_k)$ (so that $s=\sum_{k=1}^\infty
s_k$). When $k\not=n$ we have (using independence)
\[
\mathbb
E\left[\left(W_{A_k}^{(\sigma)}\right)^2\left(W_{A_n}^{(\sigma)}\right)^2\right]=s_ks_n,
\]
This is utilized so as to compute the right hand side of \eqref{new2}. We get
\[
\begin{split}
\mathbb
E\left[|\sigma(A)-\sum_{k=1}^\infty\left(W_{A_k}^{(\sigma)}\right)^2|^2\right]&=
s^2-2s^2+\sum_{k=1}^\infty \mathbb
E\left[\left(W_{A_k}^{(\sigma)}\right)^4\right]+2\sum_{k<n}s_ks_n\\
&=-s^2+3\sum_{k=1}^\infty s_k^2+2\sum_{k<n}s_ks_n\\
&=2\sum_{k=1}^\infty s_k^2,
\end{split}
\]
which is the desired right hand side of \eqref{new2}.
\end{proof}

\begin{Dn}
\label{defnew2} Let $(M,\mathscr F,\sigma)$ and
$(W_A^{(\sigma)})_{A\in\mathscr F}$ be as above. For
$A\in\mathscr F$,  we denote by ${\rm PAR}(A)$ the set of all
$\mathscr F$-partition of $A$, and we denote by
\begin{equation}
\label{new3} {\rm Var}^{(\sigma)}_{-}(A)=\inf_{(A_k)_{k\in\mathbb
N}\in{\rm PAR}(A)}\sum_{k=1}^\infty (\sigma(A_k))^2
\end{equation}
the lower variation of sum of squares.
\end{Dn}

\begin{Cy}
\label{cornew3} Let $A\in\mathscr F$, and consider the sum of
random variables squares 
\[
\sum_{k=1}^\infty \left(W_{A_k}^{(\sigma)}\right)^2
\]
on $\Omega$ for all partitions of $A$. It is $\chi^2$-distributed
on the probability space $(\Omega,\mathscr B,P)$ of Lemma
\ref{lemma3}. Then the following two conditions are equivalent.
Given $\e>0$,\\
$(i)$\hspace{1cm} $(A_k)\in{\rm PAR}(A)$ satisfies
${\rm Var}^{(\sigma)}_{-}(A)=0$ and $\sum_{k}(\sigma(A_k))^2\le \epsilon$,\\ and\\
$(ii)$\hspace{1cm}$\mathbb
E\left[|\sigma(A)-\sum_{k}\left(W_{A_k}^{(\sigma)}\right)^2|^2\right]\le
2\e$. \label{cor39}
\end{Cy}

\begin{Rk}\label{newrk}
{\rm \mbox{}\\
$(a)$ The meaning of $(ii)$ is the assertion that, when ${\rm
Var}^{(\sigma)}_{-}=0$, the random variable
\[
\sum_{k=1}^\infty \left(W_{A_k}^{(\sigma)}\right)^2
\]
is constant $P$-a.e. on $\Omega$. (It is called the quadratic
variation of $W_A^{(\sigma)}$).\\
$(b)$ The conditions in the Corollary are satisfied if
$(M,\mathscr F,\sigma)$ is taken to be the real line with the
Borel sets and the Lebesgue measure. For example if $A=[0,T]$ we
consider a sequence of partitions consisting of dyadic intervals
and
\[
{\rm Var}^{(\sigma)}_{-}([0,T])\le
T^2\sum_{k=1}^{2^n}\left(\frac{1}{2^n}\right)^2=\frac{T^2}{2^n}\rightarrow 0
\]
as $n\rightarrow\infty$.\\
$(c)$ More generally, consider $(M,\mathscr F,\sigma)$ as above,
and let $A\in\mathscr F$ be such that $0<\sigma(A)<\infty$.
Suppose that for all $n\in\mathbb N$ there is a partition
$\left\{A_1,\ldots, A_n\right\}$ such that
\begin{equation}
\label{new4} \sigma(A_k)=\frac{1}{n}\sigma(A),\quad k=1,\ldots, n.
\end{equation}
Then
\begin{equation}
\label{new5} {\rm Var}^{(\sigma)}_{-}(A)=0.
\end{equation}
Indeed, if \eqref{new4} holds, then
\[
\sum_{k=1}^n\sigma(A_k)^2=\frac{\sigma(A)^2}{n}\rightarrow 0,
\]
so \eqref{new5} follows.}
\end{Rk}

\begin{Dn}
\label{newdf1} Let $(M,\mathscr F,\sigma)$ be specified as above.
We say that it is \underline{refinable} if for every
$A\in\mathscr F$ and every $\e>0$ there exists $(A_k)\in{{\rm
PAR}}(A)$ such that
\begin{equation}
\label{newn1} |(A_k)|=\sup_{k}\sigma(A_k)<\e.
\end{equation}
\end{Dn}

\begin{Cy} Let $(M,\mathscr F,\sigma)$ and $(W_A^{(\sigma)})_{A\in\mathscr
F}$ be as above, and assume that $(M,\mathscr F,\sigma)$ is
refinable. Then for every $A\in\mathscr F$ and every
$(A_k^{(n)})\in{\rm PAR}(A)$ such that
\[
\lim_{n\rightarrow\infty}|A_k^{(n)}|=0,
\]
we have
\begin{equation}
\label{newn2} \lim_{n\rightarrow\infty}\sum_k
\left(W^{(\sigma)}_{A_k^{(n)}}\right)^2=\sigma(A),\quad P\, a.e.
\end{equation}
namely, in the limit, the left hand side random variable reduces to the constant $\sigma(A)$.
\label{cy314}
\end{Cy}
\begin{proof} Using Proposition \ref{prop37} and Corollary
\ref{cor39} we only need to show that
\begin{equation}
\label{newn3}
\lim_{n\rightarrow\infty}\sum_{k}\left(\sigma(A_k^{(n)})\right)^2=0.
\end{equation}
But the following holds for the left hand side of \eqref{newn3}:
\[
\sum_{k}\left(\sigma(A_k^{(n)})\right)^2\le|A_k^{(n)}|\sum_{k}\sigma(A_k^{(n)}).
\]
We now use that, by assumption,
\[
\lim_{n\rightarrow\infty}|A_k^{(n)}|=0,
\]
and that $\sum_{k}\sigma(A_k^{(n)})=\sigma(A)<\infty$ since
$(A_k^{(n)})\in{\rm PAR}(A)$ for all $n\in\mathbb N$. The conclusion  \eqref{newn3} now follows.
\end{proof}
\begin{Rk}{\rm It is obvious that if $(M,\mathscr F,\sigma)$ is
refinable, it has to be \underline{non-atomic}.} \end{Rk}

\begin{Dn}
\label{lalaldo}
The stochastic process $Y$ defined on $M$ with values in $\mathbf L_2(\Omega,\mathcal B,P)$ is called $\mathcal F$-adapted if the following condition
holds for every $A,B\in \mathcal F$ be such that $A\cap B=\emptyset$. Set $\mathcal F_B$ to be the sigma-algebra generated by the random variables $W^{(\sigma)}_C$, where $C$ runs through all subsets of $B$ which belong to $\mathcal F$. 
For all $x\in A$, $Y_x$ is  $\mathcal F_A$-measurable,
and, in addition, is independent of $\mathcal F_B$.
\end{Dn}

\begin{La}
\label{lalala}
Let $Y$ be a $\mathcal F$-adapted process such that the function $x\mapsto E\left(|Y(x)|^2\right)$ is measurable and $\int_M\mathbb E\left(|Y(x)|^2\right)d\sigma(x)<\infty$. Then the random variables
\[
\sum_{k=1}^nY(x_k)W^{(\sigma)}_{A_k}
\]
where $\left\{A_k\right\}_{k=1}^n$ is a covering of $M$ by pairwise disjoint measurable sets and $x_k\in A_k$, $k=1,\ldots, n$,
converges to an element in $\mathbf L_2(\Omega,\mathcal B,P)$, which we denote $\int_M Y(x)dW_x^{(\sigma)}$. Furthermore,
we have the Ito-isometry property
\begin{equation}
\label{chloe}
\mathbb E\left(\big|\int_M Y(x)dW_x^{(\sigma)}\big|^2\right)=\int_M\mathbb E\left(|Y(x)|^2\right)d\sigma(x).
\end{equation}
\end{La}

\begin{proof}
When computing the difference between two such sums, say  $\left\{A_k\right\}_{k=1}^n$ (with points $\left\{x_k\right\}_{k=1}^n$)
and $\left\{B_k\right\}_{k=1}^m$ (with points $\left\{y_k\right\}_{k=1}^m$)
we build a covering of pairwise disjoint measurable sets from the two given covering, say 
$\left\{C_k\right\}_{k=1}^p$ (with points $\left\{z_k\right\}_{k=1}^p$)
The integral
\[
\mathbb E\left(
\big|\sum_{k=1}^nY(x_k)W^{(\sigma)}_{A_k}-\sum_{k=1}^mY(y_k)W^{(\sigma)}_{B_k}\big|^2\right)
\]
can then be divided into two groups of terms: The sum
\[
\mathbb E\left(
\big|\sum_{k=1}^m(Y(z_{i_k})-Y(z_{j_k})W^{(\sigma)}_{C_k}\big|^2\right)=\sum_{k=1}^p\mathbb E\left(|Y(z_{i_k})-
Y(z_{j_k})|^2\right)\sigma(C_k)
\]
which goes to $0$ since $\int_M\left(\mathbb E\left(|Y(x)|^2\right)\right)d\sigma(x)<\infty$
and the cross-products
\[
\sum_{k,\ell=1}^p
\mathbb E\left(Y(z_k)W_{C_k}^{(\sigma)}\overline{Y(z_\ell)}W_{C_\ell}^{(\sigma)}\right)=
\big|\mathbb E\left(\sum_{k=1}^pY(z_k)W_{C_k}^{(\sigma)}\right)\big|^2
\]
which goes to $0$ since $\sigma$ is refinable, and hence non-atomic.
\end{proof}

\begin{Cy}
\label{cy322}
Let $(M,\mathscr F,\sigma)$ be a sigma-finite measure space and assume it is refinable (see Definition \ref{newdf1}). Let
$(W_A^{(\sigma)})_{A\in\mathscr F}$ be the corresponding Gaussian process (see Definition \ref{def1}). Let $f:\,\mathbb R\,\longrightarrow\,\mathbb R$ be a given $C^2$-function. Then for all $A\in\mathscr F$ such that $0<\sigma(A)<\infty$ we have:
\begin{equation}\label{itoformula}
f(W^{(\sigma)}_A)-f(0)=\int_Af^\prime(W^{(\sigma)}_x)dW_x^{(\sigma)}+\frac{1}{2}\int_Af^{\prime\prime}
(W^{(\sigma)}_x)d\sigma(x).
\end{equation}
\end{Cy}

\begin{proof}[A sketch of the proof] First note that all the terms in \eqref{itoformula} are random variables. Given $W^{(\sigma)}_A$,
by $f(W_A^{(\sigma)})$ we mean the composition of the function $W_A^{(\sigma)}$ from $\Omega$ to $\mathbb R$ with $f$, and similarly for the terms under the integrals on the right hand side of \eqref{itoformula}. Further, we stress that the first term
$\int_Af^\prime(W^{(\sigma)}_x)dW_x^{(\sigma)}$ is an Ito integral in the sense of Lemma \ref{lalala} and Remark \ref{newrk}, but now with the random process
$f^\prime(W^{(\sigma)}_x)$ occuring under the integral. By the arguments of Lemma \ref{lemma3}, we have
\[
\mathbb E_{P_\sigma}\left[\big|\int_Af^\prime(W^{(\sigma)}_x)dW_x^{(\sigma)} \big|^2\right]=\int_A\mathbb E_{P_\sigma}
\left[|f^\prime(W^{(\sigma)}_x)|^2\right]d\sigma(x).
\]
With that we note that the proof of \eqref{itoformula} is concluded through the use of the same arguments utilized in the
 proof of the classical Ito formula, see e.g. \cite{Revuz91}. Specifically, these
include $(i)$ integration by parts; $(ii)$ a stopping/truncation argument, enabling to prove the result
for $W_A$ restricted to compacts; $(iii)$ then, on compacts, any $C^2$
function may be written as a limit of
polynomial functions; and, finally, $(iv)$ the use of standard convergence
together with Proposition \ref{prop37} and 
a probabilistic Dominated Convergence allows to complete the proof.
\end{proof}

\begin{Rk}{\rm We refer to \cite[Theorem 8.2]{MR2793121} for a Ito formula where the stochastic term is computed as a Wick product integral, in the setting of an associated Gelfand triple.}
\end{Rk}

\begin{Rk}
{\rm
We note that Ito integration is done with stochastic integrands, say $Z_t$,
being adapted with respect to an underlying filtration, namely, an increasing
sequence of sub-sigma algebras ${\mathcal F_t}$. This is to say that for all $t\ge 0$, $Z_t$ is $\mathcal F_t$-measurable.
Given the fact that no natural order may be invoked in the present general setting, such terminology
obviously becomes irrelevant here. This is where Definition \ref{lalaldo} is called for.}
\end{Rk}
\subsection{Independent standard Gaussian summands}
\label{sec3.3}
\begin{Ex}
\label{ex5} Take $M=\mathbb R$ and let $\mathscr F=\mathscr
B(\mathbb R)$ denote the Borel sigma-algebra. Let $\sigma$ be a
tempered measure, that is, subject to \eqref{8},
so that we have the Gelfand triple with continuous inclusions
\begin{equation}
\label{8b} \mathcal S\hookrightarrow \mathbf
L_2(\sigma)\hookrightarrow\mathcal S^\prime
\end{equation}
where $\mathcal S$ denotes the Schwartz functions and $\mathcal
S^\prime$ denotes the tempered Schwartz distributions. Let
$\mathscr C$ be the $\mathcal S-\mathcal S^\prime$ cylinder
sigma-algebra of subsets of $\mathcal S^\prime$. We then take
\[
(\Omega,\mathscr C)=(\mathcal S^\prime, \mathscr C),
\]
and we note that the corresponding Gaussian process
$(W^{(\sigma)})$ on $(\mathcal S^\prime,\mathscr C,P)$, determined
by
\begin{equation}
\label{8c} \mathbb
E\left[e^{iW^{(\sigma)}(\varphi)}\right]=e^{-\frac{1}{2}\int_{\mathbb
R} |\widehat{\varphi}(u)|^2d\sigma(u)},\quad \varphi\in\mathcal S,
\end{equation}
satisfies the conditions in Definition \ref{def1} (and of Lemmas
\ref{lem1}-\ref{lemma3}).
\end{Ex}
\begin{Dn}
\label{def6} Let $\gamma_1$ denote the standard $\mathcal N(0,1)$ Gaussian density,
\[
d\gamma_1(x)=\frac{1}{\sqrt{2\pi}}e^{-\frac{x^2}{2}}dx\quad
\mbox{on $\mathbb R$},
\]
and set
\begin{equation}
\label{10} \Omega_\gamma=\times_{\mathbb N}\mathbb R=\mathbb
R\times \mathbb R\times\cdots,
\end{equation}
the infinite Cartesian product, with product measure
\begin{equation}
\label{11} Q\stackrel{\rm def}{=}\times_{\mathbb
N}\gamma_1=\gamma_1\times\gamma_1\times \cdots
\end{equation}
(see \cite{Hida_BM,Ka61}), defined on the cylinder sigma-algebra
$\mathscr C_\gamma$ in $\Omega_\gamma\stackrel{\rm
def}{=}\times_{\mathbb N}\mathbb R$.
\end{Dn}

Specifically, for $f_n$ a measurable
function on $\mathbb R^n$, then on $\Omega_\gamma$ set
$x=(x_1,x_2,\ldots)$ with $x_k\in\mathbb R$ ($k\in\mathbb N$),
\begin{equation}
\label{11b} F(x)=f_n(x_1,\ldots, x_n),
\end{equation}
Then, with this function $f_n$ we set:
\begin{equation}
\label{12} \begin{split}\int_{\Omega_\gamma}FdQ&=\int_{\mathbb
R^n}f_n(x_1,\ldots, x_n)d\gamma_1(x_1)\cdots d\gamma_1(x_n)\\
&=\frac{1}{(2\pi)^{n/2}}\int_{\mathbb R^n}f_n(x_1,\ldots,
x_n)e^{-\frac{1}{2}\sum_{k=1}^nx_k^2}dx_1\cdots dx_n.
\end{split}
\end{equation}

\begin{Rk}{\rm The functions $F$ in \eqref{11b} are called {\it
cylinder functions}. If the conditions in \eqref{12} hold for all
$n$, then $Q$ is uniquely determined.}
\end{Rk}

\begin{Tm}
\label{thm7} \mbox{}\\
$(i)$ Every associated Gaussian process
$(W_A^{(\sigma)})_{A\in\mathscr F}$ (see Definition \ref{def1}), corresponding to a fixed
sigma-finite measure space $(M,\mathscr F,\sigma)$ may be
realized on $(\Omega_\gamma,\mathscr C_\gamma,Q)$, i.e. in the
canonical infinite Cartesian product measure space of Definition
\ref{def6}.\\
$(ii)$ Given $(M,\mathscr F,\sigma)$ and $W^{(\sigma)}$,
$(\Omega,\mathscr C,P)$ as in Definition \ref{def1}, the
realization may be "simulated" by a system of independent,
identically distributed (i.i.d.) $\mathcal N(0,1)$ random variables
$Z_1,Z_2,\ldots$. If $(\varphi_k)_{k\in\mathbb N}$ is an
orthonormal basis of real-valued functions in $\mathbf
L_2(M,\mathscr F,\sigma)$ then a realization of $W^{(\sigma)}$ in
$(\Omega_\gamma,\mathscr C_\gamma,Q)$ may be written as
\begin{equation}
\label{13}
W_A^{(\sigma)}(\cdot)=\sum_{k=1}^\infty\left(\int_A\varphi_k(x)d\sigma(x)\right)
Z_k(\cdot),
\end{equation}
for all $A\in\mathscr F$, where $Z_k((x_i)_{i\in\mathbb
N})\stackrel{\rm def}{=}x_k$ and $(x_i)_{i\in\mathbb
N}\in\times_{\mathbb N}\mathbb R$.
\end{Tm}
\begin{proof} We first note that equation \eqref{13} is a
generalized Karhunen-Lo\`eve formula. We proceed to prove that
the expression in \eqref{13} satisfies the condition in
Definition \ref{def1}. To make sense of \eqref{13} we shall use
the Parseval formula in $\mathbf L_2(M,\mathscr F,\sigma)$. In
particular for $A\in\mathscr F$ we have
\begin{equation}
\label{14} \sigma(A)=\int_M\chi_A(x)d\sigma(x)=\sum_{k\in\mathbb
N}|\langle\chi_A,\varphi_k\rangle_\sigma|^2=\sum_{k\in\mathbb
N}|\int_A\varphi_k(x)d\sigma(x)|^2.
\end{equation}
It follows that the sum in \eqref{13} is convergent in
$\mathbf L_2(\Omega_\gamma,\mathscr C_\gamma,Q)$, and in fact is
convergent pointwise a.e. with respect to $Q$. With the use of
the dominated convergence theorem we obtain:
\[
\mathbb E_Q\left[|\sum_{k\in\mathbb
N}\int_A\varphi_k(x)d\sigma(x)Z_k(\cdot)|^2\right]=
\underbrace{\sum_{k\in\mathbb
N}|\int_A\varphi_k(x)d\sigma(x)|^2}_{\text{(by
\eqref{14})}}=\sigma(A),
\]
where we have used
\begin{equation}
\mathbb E_Q[Z_nZ_m]=\delta_{n,m}\,\quad \forall n,m\in\mathbb N.
\label{15}
\end{equation}
Hence the sum representation of $W^{(\sigma)}$ in \eqref{13} is
a well defined random variable belonging to $\mathbf
L_2(\Omega_\gamma,\mathscr C_\gamma,Q)$. For $A,B\in\mathscr F$
we have
\[
\begin{split}
\mathbb
E_Q\left[W_A^{(\sigma)}W_B^{(\sigma)}\right]&=\sum_{k\in\mathbb
N}\left(\int_A\varphi_k(x)d\sigma(x)\right)\left(\int_B\varphi_k(x)d\sigma(x)\right)\\
&=\sum_{k\in\mathbb N}\langle \chi_A,\varphi_k\rangle_\sigma
\langle \chi_B,\varphi_k\rangle_\sigma\\
&=\langle \chi_A,\chi_B\rangle_\sigma\quad (\text{by Parseval's identity})\\
&=\sigma(A\cap B),
\end{split}
\]
as required in part $(2)$ of Definition \ref{def1}. Verifying the
remaining properties in Definition \ref{def1} is immediate.
\end{proof}

\begin{Rk}{\rm
Note in particular that in representation \eqref{13} of the random
variable $W_A^{(\sigma)}$ is independent of the given choice of
an orthonormal basis in $\mathbf L_2(\sigma)$. If $(M,\mathscr
F,\sigma)=(\mathbb R,\mathscr B,\lambda)$, with $\lambda$ being
the Lebesgue measure on $\mathbb R$, then one can take as the underlying
orthonormal basis a standard wavelet, for instance the Haar
wavelets in $\mathbf L_2(\mathbb R)$. }
\end{Rk}

\begin{La}
\label{lemma8} Let $(M,\mathscr F,\sigma)$ be a sigma-finite
measure space, and let $W^{(\sigma)}$ and $(\Omega,\mathscr C,P)$
be an associated Gaussian process, see Definition \ref{def1} above.
Then for every orthonormal basis $(\varphi_k)_{k\in\mathbb N}$ in
$\mathbf L_2(M,\mathscr F,\sigma)$ the random variables
\begin{equation}
\label{16} X_k^{(\sigma)}=W^{(\sigma)}(\varphi_k),\quad
k\in\mathbb N
\end{equation}
are a collection of i.i.d. $\mathcal N(0,1)$ random variables and
\begin{equation}
\label{17}
W_A^{(\sigma)}(\cdot)=\sum_{k=1}^\infty\left(\int_A\varphi_k(x)
d\sigma(x)\right)X_k^{(\sigma)}(\cdot).
\end{equation}
\end{La}
\begin{proof} Note that the right side of \eqref{16} is the Ito
integral (see \eqref{7} in Notation \ref{n4}). Since
$W^{(\sigma)}$ is an associated Gaussian process, it follows that
$X_k^{(\sigma)}$ in \eqref{16} is Gaussian for all $k\in\mathbb
N$, with $\mathbb E_P[X_k^{(\sigma)}]=0$. Moreover, by \eqref{4} in
Lemma \ref{lemma3} we have
\begin{equation}
\label{18}
\begin{split}
\mathbb
E_P\left[X_k^{(\sigma)}X_n^{(\sigma)}\right]=\int_M\varphi_k(x)\varphi_n(x)d\sigma(x)
&=\langle \varphi_k,\varphi_n\rangle_\sigma =\delta_{m.n}
\end{split}
\end{equation}
since $(\varphi_k)_{k\in\mathbb N}$ is an orthonormal basis in
$\mathbf L_2(\sigma)$. We have already proved that the
$(X_k^{(\sigma)})_{k\in\mathbb N}$ are Gaussian, so it follows
from \eqref{18} that it is an $\mathcal N(0,1)$ i.i.d. 
collection on $(\Omega,\mathscr
C,P)$.\smallskip

Moreover, by the argument in the proof of Theorem \ref{thm7}, we may
conclude that representation \eqref{17} holds, with
convergence in $\mathbf L_2(\Omega,\mathscr C,P)$ as well as
pointwise $P$ a.e.
\end{proof}
\begin{Tm}
\label{thm9} Let $(M,\mathscr F,\sigma)$ be a sigma-finite measure
space and let $\left(W_A^{(\sigma)}\right)_{A\in\mathscr F}$ be
an associated Gaussian process on some probability
space $(\Omega,\mathscr C,P)$. Let
$\left(X_k^{(\sigma)}\right)_{k\in\mathbb N}$ be the i.i.d.
$\mathcal N(0,1)$ collection from Lemma \ref{lemma8}. Set
\[
\Gamma\,\,:\,\,\Omega\,\longrightarrow\,\,
\Omega_\gamma=\times_{\mathbb N}\mathbb R,
\]
defined by
\begin{equation}
\label{19} \Gamma(\w)=\left(X_k^{(\sigma)}(\w)\right)_{k\in\mathbb
N}\in\times_{\mathbb N}\mathbb R,\quad \w\in\Omega.
\end{equation}
Then $\Gamma$ is measurable and
\begin{equation}
P(\Gamma^{-1}(C))=Q(C) \label{20}
\end{equation}
holds for all $C\in\mathscr C_\gamma$, i.e. the sigma-algebra of
subsets in $\times_{\mathbb N}\mathbb R$ generated by the
cylinder sets.
\end{Tm}

\begin{Rk}{\rm
For \eqref{20}, we will use the notation
\begin{equation}
\label{21} P\circ\Gamma^{-1}=Q,
\end{equation}
with
\[
\Gamma^{-1}(C)=\left\{\w\in\Omega\,|\, \Gamma(\w)\in C\right\}. \]
}
\end{Rk}
\begin{proof}[Proof of Theorem \ref{thm9}]
Fix $m\in\mathbb N$, and let $I_k=(a_k,b_k]$, $1\le k\le m$ be a series of
finite intervals. It follows from \eqref{19} that if
\[
C=C_m=\left\{(x_k)_{k\in\mathbb N}\,|\, x_k\in I_k,\, 1\le k\le
m\right\}
\]
then
\[
\Gamma^{-1}(C)=\left\{\w\in\Omega\,|\, a_k\le
X^{(\sigma)}_k(\w)<b_k\right\}
\]
and so $\Gamma$, defined by \eqref{19}, is measurable. But we know
that $(X^{(\sigma)}_k)_{k\in\mathbb N}$ is an $\mathcal N(0,1)$ i.i.d. 
collection in $(\Omega,\mathscr C,P)$, and $(Z_k)_{k\in\mathbb N}$ is
an $\mathcal N(0,1)$ i.i.d. collection in $(\Omega_\gamma,\mathscr C_\gamma,
Q)$. Hence
\begin{equation}
\label{22} P(\Gamma^{-1}(C))=Q(C_m)=\prod_{k=1}^m\gamma_1(I_k).
\end{equation}
Thus the desired conclusion \eqref{20} holds for cylinders in the
respective probability spaces. On the other hand, the subset of
all sets $C\subset\Omega_\gamma$ for which \eqref{20} holds is a
sigma-algebra. It follows that \eqref{20} is valid on the
respective sigma-algebras generated by the cylinder sets.
\end{proof}
\subsection{Generalized Fourier transform}
\begin{Cy}
\label{cor10} Let $(M,\mathscr F,\sigma)$ be a sigma-finite
measure space, and let $\left(W_A^{(\sigma)}\right)_{A\in\mathscr
F}$ be an associated Gaussian process on the probability space
$(\Omega,\mathscr B, P)$. Then the closed linear span in
$\mathbf L_2(\Omega,\mathscr B, P)$ of the functions
\[
\left\{e^{iW_A^{(\sigma)}(\cdot)}\, |\, A\in\mathscr F\right\}
\]
is $\mathbf L_2(\Omega,\mathscr B, P)$.
\end{Cy}

\begin{proof}
By Theorem \ref{thm9} we may assume that
$\left(W_A^{(\sigma)}\right)_{A\in\mathscr F}$ is realized in
$\mathbf L_2(\times_{\mathbb N}\mathbb R,\mathscr C_\gamma, Q)$,
with $Q=\times_{\mathbb N}\gamma_1$ being the infinite dimensional
measure on $\Omega_\gamma=\times_{\mathbb N}\mathbb R$; see also
Theorem \ref{thm7}. But in this infinite-product space the
conclusion is clear. If $F\in\mathbf L_2(\times_{\mathbb
N}\mathbb R,\mathscr C_\gamma, Q)$ satisfies
\[
\int_{\times_{\mathbb N}\mathbb R} F e^{iW_A}dQ=0,\quad\forall
A\in\mathscr F,
\]
then a direct computation, using \eqref{13}, shows that $F=0$
a.e. $Q$.
\end{proof}

\begin{Dn}
\label{def11} Let $(M,\mathscr F,\sigma)$ be a sigma-finite
measure space, and consider the positive definite function
\begin{equation}
\label{23}
K^{(\sigma)}(A,B)=e^{-\frac{\sigma(A)+\sigma(B)}{2}}e^{\sigma(A\cap
B)}=e^{-\frac{1}{2}\|\chi_A-\chi_B\|^2_\sigma},\quad
A,B\in\mathscr F.
\end{equation}
We shall denote the corresponding reproducing kernel Hilbert
space by $\mathcal H^{(\sigma)}$.
\end{Dn}
\begin{Cy}
\label{cor12} (The generalized Fourier transform.) Let
$\left(W_A^{(\sigma)}\right)_{A\in\mathscr F}$ be a Gaussian
process associated with the sigma-finite measure space  $(M,\mathscr
F,\sigma)$, and let $(\Omega,\mathscr B, P)$ be the underlying
probability space. Let $\mathcal H^{(\sigma)}$ be the reproducing
kernel Hilbert space as in Definition \ref{def11}. For
$F\in\mathbf L_2(\Omega,\mathscr B, P)$, set
\begin{equation}
\label{24} \widehat{F}(A)=\mathbb E_P\left[Fe^{iW_A}\right],\quad
A\in\mathscr F.
\end{equation}
Then the map $F\mapsto \widehat{F}$ is unitary from $\mathbf
L_2(\Omega,\mathscr B, P)$ onto $\mathcal H^{(\sigma)}$.
\end{Cy}
\begin{proof}
Using Corollary \ref{cor10}, it is enough to prove that for all
$A,B\in\mathscr F$ we have
\begin{equation}
\label{25} \mathbb
E\left[e^{-iW_B}e^{iW_A}\right]=K^{(\sigma)}(A,B)=e^{-\left\{\frac{\sigma(A)+\sigma(B)}{2}
+\sigma(A\cap B)\right\}}
\end{equation}
(see \eqref{23}). A direct computation yields
\[
\begin{split}
E\left[e^{-iW_B}e^{iW_A}\right]&=E\left[e^{i(W_A-W_B)}\right]\\
&=e^{-\frac{1}{2}\|\chi_A-\chi_B\|^2_\sigma}\quad (\text{by
Theorem \ref{thm9}})\\
&=K^{(\sigma)}(A,B)\quad\,\,\,\,\,(\text{by \eqref{23}}),
\end{split}
\]
for all $A,B\in\mathscr F$.
\end{proof}

\begin{Cy}
\label{cor13} Let $\left(W_A^{(\sigma)}\right)_{A\in\mathscr
F}$ be a collection of Gaussian random variables associated with a sigma-finite
measure space $(M,\mathscr F,\sigma)$. Let
$(\varphi_k)_{k\in\mathbb N}$ be an orthonormal basis of $\mathbf
L_2(M,\mathscr F,\sigma)$ (we denote this space by $\mathbf
L_2(\sigma)$), and let
$
\Gamma\,\,:\,\,\Omega\,\,\longrightarrow\,\,\times_{\mathbb
N}\mathbb R\stackrel{\rm def.}{=}\Omega_\gamma$, with
\begin{equation}
\label{26} \Gamma(
\w)=(X_k^{(\sigma)}(\w))_{k\in\mathbb N}.
\end{equation}
being the corresponding coordinate system, i.e.
\[
X_k^{(\sigma)}\stackrel{\rm
def.}{=}W^{(\sigma)}(\varphi_k)=\int_M\varphi_k(x)dW^{(\sigma)}(x)
\]
being the Ito-integral representation of \eqref{7}. Then
\begin{equation}
\label{27} \Gamma\left(\mathbf
L_2(\sigma)\right)\subset\ell^2\subset\times_{\mathbb N}\mathbb R
\end{equation}
\end{Cy}

\begin{proof} Recall that we consider real-valued functions. Let
$\langle\cdot, \cdot\rangle_\sigma$ denote the standard inner
product in $\mathbf L_2(\sigma)$. Then for all $f\in \mathbf
L_2(\sigma)$ we have
\[
W^{(\sigma)}(f)\in\mathbf L_2(\Omega, P),
\]
and, by Lemma \ref{lemma3},
\[
\Gamma(f)=\left(\langle
f,\varphi_k\rangle_\sigma\right)_{k\in\mathbb
N}\in\ell_2\subset\times_{\mathbb N}\mathbb R,
\]
and
\begin{equation}
\label{28} \mathbb E_P\left[\left(W^{(\sigma)}(f)\right)^2\right]=\|f\|_\sigma^2=
\sum_{k\in\mathbb N}|\langle f,\varphi_k\rangle|^2,
\end{equation}
and
\begin{equation}
\label{29} W^{(\sigma)}(f)=\sum_{k\in\mathbb N}\langle
f,\varphi_k\rangle_kX_k^{(\sigma)}
\end{equation}
is well defined, Gaussian, and satisfies \eqref{28}.
\end{proof}
\begin{Tm}
\label{thm14} Let $(M,\mathscr F,\sigma)$,
$\left(W_A^{(\sigma)}\right)_{A\in\mathscr F}$ and
$(\Omega,\mathscr B, P)$ be as in Corollary \ref{cor13}. Then the
measure $P$ is quasi-invariant with respect to the $\mathbf
L_2(\sigma)$ translations in $\Omega$. More precisely, we have
\begin{equation}
\label{30}
\int_\Omega  F(\cdot+f)dP(\cdot)=\int_\Omega F(\cdot)e^{-\frac{1}{2}\|f\|_\sigma^2+W^{(\sigma)}(f)(\cdot)}dP(\cdot)
\end{equation}
for all $F\in\mathbf L_2(\Omega, P)$ and all $f\in\mathbf L_2(\sigma)$.
\end{Tm}

\begin{proof}
Note that in the formulation of \eqref{30} we make use of the coordinate system
\[
\Gamma\,:\, \Omega\,\,\longrightarrow\,\, \times_{\mathbb N}\mathbb R
\]
from Theorem \ref{thm9} and Corollary \ref{cor13}.  Hence the term
$e^{-\frac{1}{2}\|f\|_\sigma^2+W^{(\sigma)}(f)(\cdot)}$
is the associated Radon-Nikodym derivative. Using Theorem \ref{thm9}, and passing to cylinder functions, we note that
\eqref{30} follows from a computation  of $\mathbb R^m$-integrals
for all $m\in\mathbb N$. We check that if $F_m$ is $\mathbf L_2$ with respect to the standard $\mathbb R^m$-Gaussian joint density
\[
\gamma_m=\underbrace{\gamma_1\times \cdots\times\gamma_1}_{\text{$m$ times}},
\]
then
\begin{equation}
\label{31}
\begin{split}
\iint_{\mathbb R^m}F_m(x_1+\langle f,\varphi_1\rangle_\sigma,\ldots,
x_m+\langle f,\varphi_m\rangle_\sigma)
d\gamma_m(x_1,\ldots, x_m)&=\\
&\hspace{-6.5cm}= \iint_{\mathbb
R^m}F_m(x_1,\ldots, x_m)e^{-\frac{1}{2}\sum_{k=1}^m\langle
f,\varphi_k\rangle^2_\sigma}e^{\sum_{k=1}^mx_k\langle f, \varphi_k
\rangle_\sigma}d\gamma_m(x_1,\ldots, x_m).
\end{split}
\end{equation}
From Theorem \ref{thm9} and Corollary \ref{cor13} we have that
\begin{equation}
\label{32}
W^{(\sigma)}(f)(\cdot)=\sum_{k=1}^\infty \langle f,\varphi_k\rangle_\sigma X_k^{(\sigma)}(\cdot)
\end{equation}
and
\begin{equation}
\label{33}
\mathbb E_\sigma\left[\left(W^{(\sigma)}(f)\right)^2\right]
=\sum_{k=1}^m\langle f,\varphi_k\rangle_\sigma^2=\|f\|^2_{\mathbf L_2(\sigma)}.
\end{equation}
Using again Theorem \ref{thm9} and Kolmogorov's induction limit construction, we finally note that the desired formula \eqref{30}
follows from \eqref{31}-\eqref{33} above.
\end{proof}

Our quadratic variation result in the first part of the present
section  is motivated by, and is a generalization of a
classical theorem for the Brownian motion, often called L\'evy's
theorem, see e.g. \cite{Buh62}. Similarly, our decomposition
theorem (Theorem \ref{thm7}) is motivated by, and extends a
classical result often called  a {\it Kahrunen-Lo\`eve
decomposition}, see e.g. \cite{Dur98}. In the above,
we have made use of the theory of Ito-integration and its generalizations; see e.g. \cite{MR0388538,MR3051797,MR2508233,MR2454534,MR2444857,MR2397796}
\section{Gaussian stochastic calculus}
\label{sec4}
In this section, a preparation for Section \ref{sec6},  we return to the Gaussian process of Section \ref{sec2},  
applying results of Section \ref{sec3}.
 
\setcounter{equation}{0}
\subsection{Gaussian processes and tempered measures: The Minlos theorem.}

Let $\sigma$ be a tempered measure. The map
\begin{equation}
\label{nnew5}
\varphi\,\mapsto\,e^{-\frac{1}{2}\int_{\mathbb R}|\widehat{\varphi}(u)|^2d\sigma(u)}
\end{equation}
is continuous and positive definite on $\mathcal S$ (see \cite[Proposition 3.3, p. 714]{ajnfao}). 
An application of Minlos' theorem (see \cite{MR0435834}) yields:
\begin{Tm}
For every tempered measure $\sigma$ there exists a uniquely defined probability measure $P_\sigma$ on $\mathcal S^\prime$, defined on the
sigma-algebra $\mathscr C=\mathscr C(\mathcal S^\prime)$ generated by the cylinder sets (see \eqref{eqdefcyl}), and determined
by the following condition:
\begin{equation}
\label{new6}
\mathbb E_{P_\sigma}\left[e^{iX_\varphi^{(\sigma)}}\right]=
\int_{\mathcal S^\prime}e^{i\langle \varphi,\xi\rangle}dP_\sigma(\xi)=e^{-\frac{1}{2}\int_{\mathbb R}|\widehat{\varphi}(u)|^2d\sigma(u)}.
\end{equation}
\end{Tm}
Using the corresponding expectation
\[
\mathbb E_\sigma[U]=\int_{\mathcal S^\prime} U(\w)P_\sigma(\w),
\]
\eqref{new6} takes the equivalent form
\[
\mathbb E_\sigma\left[e^{i\langle \cdot,\varphi\rangle}\right]=e^{-\frac{1}{2}\|\widehat{\varphi}\|^2_{\sigma}}.
\]
In the above we make use of reproducing kernel Hilbert spaces and their association with the study of Gaussian processes; 
see e.g.
\cite{MR2952544,MR2805533,MR2788708,MR2599216,MR2459226}.
\subsection{A Gaussian process realized on $\mathcal S^\prime$.}
\begin{La}
The pairing \eqref{pairing} between $\mathcal S$ and $\mathcal S^\prime$ extends to a pairing between
$\mathbf L_2(\mathbb R)$ and $\mathcal S^\prime$. The stochastic process
\[
X^{(\sigma)}_\varphi(\xi)=\langle\xi,\varphi\rangle,\quad \varphi\in\mathcal S,
\]
may be extended to a process
\[
X^{(\sigma)}_f(\xi)=\langle\xi,f\rangle,\quad f\in\mathbf L_2(\mathbb R,dx).
\]
\label{la41}
\end{La}

\begin{proof}
Set
\begin{equation}
\label{old7}
X_t^{(\sigma)}(\xi)=\langle\xi,\chi_{[0,t]}\rangle.
\end{equation}
$X_t^{(\sigma)}$ is a zero mean Gaussian process with variance
\begin{equation}
\label{realform}
\int_{\mathbb R}|\widehat{\chi_{[0,t]}}(u)|^2d\sigma(u)=2\int_{\mathbb R}\frac{1-\cos(ut)}{u^2}d\sigma(u)=4
\int_{\mathbb R}\frac{\sin^2(\frac{ut}{2})}{u^2}d\sigma(u).
\end{equation}
\end{proof}

In the following lemma we compute the covariance and related quantities. In formulas \eqref{f1} and \eqref{f2} 
below we set
\begin{equation}
\label{f0}
r_\sigma(t)=\mathbb E_\sigma\left[|X_t^{(\sigma)}|^2\right].
\end{equation}

\begin{La} We set $X_0^{(\sigma)}=0$.
The following formulas hold:\\
\begin{eqnarray}
\mathbb E_\sigma\left[|X_t^{(\sigma)}|^2\right]&=&4\int_{\mathbb R}\frac{\sin^2\left(\frac{ut}{2}\right)}{u^2}d\sigma(u),\\
\label{f1}
\mathbb E_\sigma\left[|X_t^{(\sigma)}-X_s^{(\sigma)}|^2\right]&=&r_\sigma(s-t),\\
\label{f2}
\mathbb E_\sigma\left[X_t^{(\sigma)}X_s^{(\sigma)}\right]&=&\frac{r_\sigma(t)+r_\sigma(s)-r_\sigma(s-t)}{2}.
\end{eqnarray}
\label{la42}
\end{La}
The computations are standard and will be omitted.\\

Thus $X^{(\sigma)}$ is a stationary-increment Gaussian process. Recall that these are (not necessarily real) 
Gaussian processes with covariance function of the form
\begin{equation}
\label{sincre}
2\mathbb E_\sigma\left[X_t^{(\sigma)}X_s^{(\sigma)}\right]=r(t)+\overline{r(s)}-r(t-s),\quad t,s\in\mathbb R,
\end{equation}
Functions $r$ for which the kernel \eqref{sincre} is positive definite on the real line have been investigated by Schoenberg, von Neumann 
and Krein; see \cite{MR0004644,MR0012176}. They are of the form
\[
r(t)=r_0+ict-\int_{\mathbb R}\left\{e^{itv}-1-\frac{itv}{v^2+1}\right\}d\sigma(v)
\]
where $r_0\ge 0$, $c\in\mathbb R$ and $\sigma$ is a tempered measure satisfying \eqref{8} for $p=1$. 
In the case where $\sigma$ is even, one recovers formula \eqref{realform}.\smallskip

We note that such processes admit derivatives which are generalized stationary processes
with covariance $\widehat{\sigma}(t-s)$, with the Fourier transform being computed in the sense of distributions. See 
\cite{aal2,MR2793121,ajnfao} for more information.

\section{A Gaussian process realized on $C(\mathbb R)$.}
\setcounter{equation}{0}
\label{sec6}
We are now ready to introduce the two realizations of probability spaces for the family of Gaussian processes considered above.
\begin{Tm}
\label{tm61}
There is a unique measure $Q_\sigma$ defined on the cylinder sigma-algebra of $\Omega=C(\mathbb R)$ such that
\begin{equation}
\label{new44}
X_t^{(\sigma)}(\w)=\w(t),\quad t\in\mathbb R,\quad \w\in\Omega=C(\mathbb R),
\end{equation}
and $X^{(\sigma)}_0=0$ is a Gaussian process with zero mean and covariance
\begin{equation}
\mathbb E_{Q_\sigma}\left[X_t^{(\sigma)}X_s^{(\sigma)}\right]=\frac{r_\sigma(t)+r_\sigma(s)-r_\sigma(s-t)}{2}
\end{equation}
(see \eqref{f2}), with $r_\sigma$ as in \eqref{f0}.
\end{Tm}
\begin{proof}
For a finite set $t_1,\ldots, t_n\in\mathbb R$ of sample points and a measurable function $f_n$ on $\mathbb R^n$ we set
\begin{equation}
F(\w)=f_n(\w(t_1),\ldots, \w(t_n)),\quad \w\in\mathbb C(R).
\label{new17}
\end{equation}
For  $a>0$ define
\begin{equation}
\label{ga}
g_a(x)=\frac{1}{a\sqrt{2\pi}}e^{-\frac{x^2}{2a^2}},\quad x\in\mathbb R.
\end{equation}
and for $0<t_1<t_2<\ldots<t_n$ let
\begin{equation}
\begin{split}
\mathcal L_n(F)&=\int_{\mathbb R^n}f_n(x_1,\ldots,
x_n)g_{r_\sigma(t_1)}(x_1)g_{r_\sigma(t_2-t_1)}(x_2-x_1)\cdots
g_{r_\sigma(t_n-t_{n-1})}(x_n-x_{n-1})dx_1dx_2\cdots dx_n\\
&=\int_{\mathbb R^n}f_n(x_1,x_2+x_1,\ldots, x_n+x_{n-1}+\cdots+x_1)
g_{r_\sigma(t_1)}(x_1)g_{r_\sigma(t_2-t_1)}(x_2)\cdots
g_{r_\sigma(t_n-t_{n-1})}(x_n)dx_1dx_2\cdots dx_n,
\end{split}
\end{equation}
where $f_n$ is the function in $\mathbb R^n$ introduced in
\eqref{new17}.\smallskip

Using Kolmogorov's theorem we see that there is a unique measure
$Q_\sigma$ with the property stated in the theorem such that
for a cylinder function $F$ as in \eqref{new17} we have
\begin{equation}
\label{new19}
\int_\Omega FdQ_\sigma=\mathcal L_n(F).
\end{equation}
The other claims are easily verified. We omit the details.
\end{proof}
In this section we make use of the two general procedures for constructing Gaussian processes, and the corresponding 
probability spaces. To simplify matters, we divide them into $(i)$ the inductive limit construction first proposed by 
Kolmogorov, and $(ii)$ the
alternative approach based on Gelfand triples. For references to the first, see \cite{MR0032961,MR1805833},
for the second, see e.g., \cite{MR1810995}.\\

We further mention the following result from \cite{ajnfao} and
\cite{MR2793121} (we refer the reader to these references for the proof).

\begin{Tm}
\label{place-voltaire}
Let $\sigma$ be a tempered measure (see Definition \ref{def21}) and let $\left(X_\varphi^{(\sigma)}\right)_{\varphi\in\mathcal S}$
be the corresponding Gaussian process indexed by the real Schwartz space $\mathcal S$. Further, let $P_W$ denote the standard Wiener-measure. For $f\in\mathbf L_2(\mathbb R)$, let $\widetilde{f}$ denote the corresponding Gaussian process with zero mean and covariance
\[
\mathbb E_{P_W}(|\widetilde{f}|^2)=\|f\|_2^2.
\]
Then, there exists a continuous linear operator $Q^{(\sigma)}$ from $\mathcal S$ into $\mathbf L_2(\mathbb R)$ such that
\[
\mathbb E_{P_W}\left[e^{i\widetilde{Q^{(\sigma)}(\varphi)}}\right]=\mathbb E_{P_\sigma}\left[e^{iX_\varphi^{(\sigma)}}\right]
=e^{-\frac{1}{2}\int_{\mathbb R}|\widehat{\varphi}(u)|^2d\sigma(u)}.
\]
\end{Tm}

\section{Equivalence of measures}
\label{sec7}
In this section we provide details (Theorem \ref{thm82} below) associated with 
the equivalence of the two realizations of probability 
spaces for the family of Gaussian processes considered above.
\setcounter{equation}{0}
\subsection{Two probability spaces}
In the approach based on Gelfand triples
we realize the family of stationary-increment processes on the
probability space $(\mathcal S^\prime, \mathscr C(\mathcal
S^\prime),P_\sigma)$, while in the approach based instead on
Kolmogorov consistency theorem, a probability space is
$(C(\mathbb R), \mathscr C,Q_\sigma)$, where $\mathscr C$ is the
corresponding cylinder sigma-algebra. In both cases the Gaussian
process is indexed by a fixed tempered measure $\sigma$.\smallskip

The setting in the present section is as follows. Fix a tempered
measure $\sigma$ on $\mathbb R$ and consider the above two
variations of the $\sigma$-Gaussian process, For $\w\in C(\mathbb
R)$ let $\w^\prime$ be the corresponding Schwartz tempered
distribution.\smallskip

We shall write $X_{\rm Kolm}^{(\sigma)}(t)$ for the process in the
Kolmogorov realization, and $X_{\rm Gel}^{(\sigma)}$ for the
Gelfand-triple realization, i.e.
\begin{eqnarray}
\label{new20} X_{\rm Kolm}^{(\sigma)}(t)(\w)&=&\w(t),\quad\forall
\w\in
C(\mathbb R),\, \forall t\in\mathbb R,\\
X_{\rm
Gel}^{(\sigma)}(t)(\xi)&=&\langle\xi,\chi_{[0,t]}\rangle,\quad\forall
\xi\in\mathcal S^\prime,\, \forall t\in\mathbb R. \label{new21}
\end{eqnarray}
\subsection{Equivalence}
In the following statement recall that
cylinder sets have been defined in \eqref{eqdefcyl}.

\begin{La}
\label{lanew} Let $\sigma$ be a tempered measure as introduced in
Section \ref{sec4} (see Lemmas \ref{la41} and \ref{la42}). Then
the Gaussian process $(X_\varphi^{(\sigma)})_{\varphi\in\mathcal
S}$ determined by
\[
\mathbb E\left[X_\varphi^{(\sigma)}\right]=0\quad\text{and}\quad \mathbb
E\left[e^{iX_\varphi^{(\sigma)}}\right]=e^{-\frac{1}{2}\int_{\mathbb
R} |\widehat{\varphi}(u)|^2d\sigma(u)} \\
\]
extends from $\mathcal S$ to $\mathbf L_2(\mathbb R)$, i.e. the
$\mathcal S$-$\mathcal S^\prime$ pairing
$X_\varphi^{(\sigma)}(\xi)=\langle \xi,\varphi\rangle$ extends
uniquely to $X_f^{(\sigma)}$.
By abuse of notation we shall write
\[
X_f^{(\sigma)}(\xi)=\langle \xi,f\rangle.
\]
\end{La}
\begin{proof}
For $f\in\mathbf L_2(\mathbb R)$, interpret $X_f^{(\sigma)}(\xi)$
as the Ito-integral
\begin{equation}
\label{new456} X_f^{(\sigma)}(\xi)=\int_{\mathbb
R}f(x)dX^{(\sigma)}(x),
\end{equation}
where the right hand side of \eqref{new456} is computed with the
use of Lemma \ref{lemma3}.
\end{proof}

Fix $k\in\mathbb N$, and consider a Gaussian measure $Q$ on $C(\mathbb R^k)$ constructed from Gaussian transition probabilities, see e.g. 
\eqref{ga}-\eqref{new17} above. We denote by
\[
C_{\rm temp}(\mathbb R^k)=C(\mathbb R^k)\cap \mathcal S^\prime_{\mathbb R^k}
\]
the set of tempered continuous functions on $\mathbb R^k$ (see \cite{Treves67}). Furthermore, let 
$\stackrel{\bullet}{\mathbb R^k}=\left(\mathbb R^k\cup\left\{\infty\right\}\right)^\sim$ denote the one-point 
compactification of $\mathbb R^k$. Then we have:

\begin{La}
\label{Lanew}
It holds that
\begin{equation}
\label{eq300814}
C(\mathbb R^k)\cap\prod_{\mathbb R^k}\stackrel{\bullet}{\mathbb R^k}\subset C_{\rm temp}(\mathbb R^k).
\end{equation}
\end{La}
\begin{proof}
The functions $\w$ in the intersection on the left hand side of \eqref{eq300814} have well defined limits at $\infty$. There are in
particular bounded and hence tempered.
\end{proof}

\begin{La}
The Gaussian measure $Q_\sigma$ in Theorem \ref{tm61} is supported in $C(\mathbb R^k)\cap\prod_{\mathbb R^k}\stackrel{\bullet}{\mathbb R^k}$
\end{La}

\begin{proof}
The fact that $Q_\sigma(C(\mathbb R^k)\cap\prod_{\mathbb R^k}\stackrel{\bullet}{\mathbb R^k})=1$ is contained in 
\cite[Appendix A]{MR0161189}.
\end{proof}
\begin{Tm}
\label{thm82}
For a fixed tempered measure $\sigma$ let $X^{(\sigma)}_{\rm
Kolm}$ and $X^{(\sigma)}_{\rm Gel}$ be the associated Gaussian
processes arising from Kolmogorov and Gelfand constructions
respectively. Let
\[
T\,:\, C_{\rm temp}(\mathbb R)\,\longrightarrow\, \mathcal S^\prime
\]
be defined by
\begin{equation}
\label{new22} T(\w)=\w^\prime,\quad \w\in C_{\rm temp}(\mathbb R);
\end{equation}
(the prime in $\w^\prime$ denotes the derivative in the sense of
distributions). Then,
\begin{equation}
\label{new222}
X^{(\sigma)}_{\rm Gel}(\w^\prime)=X^{(\sigma)}_{\rm Kolm}(\w)=\w(t),
\end{equation}
and
\[
Q_\sigma\circ T^{-1}=P_\sigma,
\]
meaning that for all cylinder sets $\Delta\subset\mathcal S^\prime$, with
\[T^{-1}(\Delta)=\left\{\w\in C_{\rm temp}(\mathbb R)\,\,\text{such that}\,\, \w^\prime\in\Delta\right\},
\]
then
\begin{equation}
Q_\sigma\left(T^{-1}(\Delta)\right)=P_\sigma(\Delta).
\label{new23}
\end{equation}
\end{Tm}
\begin{proof}
To prove \eqref{new222} let $\w\in C_{\rm temp}(\mathbb R)$ and $t\in\mathbb R$. Since the process is assumed to be zero 
at $t=0$ we can write
\[
\begin{split}
X^{(\sigma)}_{\rm Gel}(\w^\prime)&=\langle\w^\prime, \chi_{[0,t]}\rangle\\
&=-\langle \chi_{[0,t]}^\prime, \w\rangle\\
&=\w(t)\\
&=X^{(\sigma)}_{\rm Kolm}(t)(\w).
\end{split}
\]

We now turn to the proof of \eqref{new23}. Note that this yields an explicit formula for the measure $P_\sigma$
on $\mathcal S^\prime$ which was previously obtained indirectly as an application of Minlos' theorem to the
right hand side of \eqref{6}, i.e. the Minlos-Gelfand approach does not yield a construction of $P_\sigma$, only existence.
Let $A\subset \mathbb R^n$ be a Borel set. Consider a cylinder set of the form \eqref{eqdefcyl}, with
\[
\varphi_k(u)=\chi_{[0,t_k]}(u),\quad k=1,\ldots, n.
\]
Then,
\[
\begin{split}
P_\sigma\left\{\w^\prime\in\mathcal S^\prime\,|\, (\langle\w^\prime,\chi_{[0,t_1]}\rangle,\ldots,
\langle\w^\prime,\chi_{[0,t_n]}\rangle)\in A\right\}=Q_\sigma\left\{\w\in C_{\rm temp}(\mathbb R)\,|\,(\w(t_1),\ldots, \w(t_n))\in
A\right\},
\end{split}
\]
where $P_\sigma$ and $Q_\sigma$ refer to the functional measures on $\mathcal S^\prime$ and $C(\mathbb R)$
respectively.\smallskip

\end{proof}

{\bf Acknowledgment:} It is a pleasure to thank the referee for his comments on the first version of the paper.
\bibliographystyle{plain}
\def\cfgrv#1{\ifmmode\setbox7\hbox{$\accent"5E#1$}\else
  \setbox7\hbox{\accent"5E#1}\penalty 10000\relax\fi\raise 1\ht7
  \hbox{\lower1.05ex\hbox to 1\wd7{\hss\accent"12\hss}}\penalty 10000
  \hskip-1\wd7\penalty 10000\box7} \def\cprime{$'$} \def\cprime{$'$}
  \def\cprime{$'$} \def\lfhook#1{\setbox0=\hbox{#1}{\ooalign{\hidewidth
  \lower1.5ex\hbox{'}\hidewidth\crcr\unhbox0}}} \def\cprime{$'$}
  \def\cprime{$'$} \def\cprime{$'$} \def\cprime{$'$} \def\cprime{$'$}
  \def\cprime{$'$}

\end{document}